\documentclass[12pt,a4paper]{article}
\usepackage[T2A]{fontenc}
\font\gothic = eufm10

\newtheorem{thm}{Theorem}
\newtheorem{lem}{Lemma}
\newtheorem{cor}{Corollary}

\newtheorem{rem}{Remark}

\newtheorem{exm}{Example}

\begin{document}

\title {Limit theorem for  randomly indexed sequence of random processes}

\author{E.  E. Permyakova}





\maketitle

{\bf Abstract}
In this paper  is proved the limit theorem for randomly indexed sequence of random processes in the case where sequences of random index and random processes are independent, also the estimation of convergence rate is obtained.

\section{Introduction}

The study of asymptotic behavior of randomly indexed sequence enjoy considerable attention in connection with applications in the  theory of queues, Markov processes, word space modeling. In the papers [1-4] the conditions of convergence and asymptotic behavior of randomly indexed random variables is studied. In this work we consider the sequences of randomly indexed random processes, defined in Skorokhod space $D[0,1]$. The obtained limit theorem makes it easy to get the limit theorems for some random processes with random substitution defined in $D[0,1]$. The convergence conditions are quite weak,  which allows to apply the theorem in many applications. The estimation of convergence rate is also obtained.

\section{Main results}
Recall the definition of the Skorokhod space and the metric in it.
Denote by
$\Delta[0,1]$  class of strictly increasing continuous mapping of the segment  $[0,1]$ on it self such that
$\lambda(0)=0, \ \lambda(1)=1$; $D[0,1]$ is a Skorokhod space, i.e. the space of functions defined in the segment $[0,1]$ and taking values in ${\rm I\kern-.23em R}$, right-continuous and with a finite limit on the left. In the space $D[0,1]$ we will consider the Skorokhod's metric $$\rho(x,y)=$$
$$\inf\{\varepsilon>0: \exists \lambda\in \Delta[0,1], \
\sup_{0\le t\le 1}|x(t)-y(\lambda(t))|\le \varepsilon,\ \
\sup_{0\le t,s\le 1, s\not=t}\left|ln\frac{\lambda(t)-\lambda(s)}{t-s}\right|\le \varepsilon\},$$ $x,y\in
D[0,1]$. It is well know \cite{Billingsley} that this metric turns $D[0,1]$ into Polish space. Then (see, for example, \cite{Bulinski}) cylindrical $\sigma$- algebra coincides with Borel $\sigma$- algebra, and the Borel $\sigma$- algebra \cite{Billingsley} is generated by maps of type $f\rightarrow f(a)$ for $a \in[0,1]$.
\begin{thm} \label{maintheorem}
Let $Y_n$ be a sequence of random processes in $D[0,1]$ such that $$Y_n \stackrel{d}{\to} Y \mbox{ as } n\to\infty,$$ and trajectories of  $Y$ are continuous, $\nu_n(t)$ is a sequence of random processes in $D[0,1]$ with non-decreasing non-negative trajectories, for every $t\in [0,1]$ $\nu_n(t)$  	
take values in  ${\rm I\kern-.23em N}$, $f(n)$ is a function taking values in ${\rm I\kern-.23em R}$ such that $f(n)\to\infty$ as $n\to\infty$ and
$$ \frac{\nu_n}{f(n)}\stackrel{d}{\to} \nu \mbox{ in } D[0,1] \mbox{ as } n\to\infty,$$ and besides it exists $c>0$ such that  $\nu(t)>c$ for all $t \in [0,1]$.  Let $Y_n$ and $\nu_n$ are independent.
Then
\begin{equation}       \label{theorem1}
Y_{\nu_n}\stackrel{d}{\to} Y \mbox{ as } n\to\infty 
\mbox{ in } D[0,1],
\end{equation}
where the convergence holds in the uniform norm.
\end{thm}
To prove this theorem we need the following technical preliminary  result.
\begin{lem}           \label{lemma}
Let $x_n(t)$ be a sequence of functions in $D[0,1]$ such that $x_n\to x$ in $D[0,1]$ as $n\to\infty$, $x$ is continuous and the sequence of c\`adl\`ag non-decreasing non-negative functions $\mu_n(t)$ for every $t\in[0,1]$  take values in ${\rm I\kern-.23em N}$ and
$$ \frac{\mu_n}{f(n)}\stackrel{d}{\to} \mu \mbox{ in } D[0,1] \mbox{ as } n\to\infty$$ for some function $f(n)$, where $f(n)\to\infty$ as $n\to\infty$ and besides it exists $c>0$ such that  $\mu(t)>c$ for all $t \in [0,1]$. Then $$x_{\mu_n}\to x \mbox{ in } D[0,1] \mbox{ as } n \to\infty,$$ where the convergence holds in the uniform norm.
\end{lem}
{\it Proof.}
Let $0<\varepsilon<c$ is arbitrary. The convergence $x_n\to x$ implies the existence of $n_1$ such that for all $n>n_1$ $\rho_D(x_n,x)<\varepsilon$.

The convergence $ \frac{\mu_n}{f(n)}\stackrel{d}{\to} \mu \mbox{ in } D[0,1] \mbox{ as } n\to\infty$ implies that exists $n_2$ such that for all $n>n_2$ it exists $\lambda_n \in \Delta[0,1]$ and $$ \sup_{0\le t\le 1}\left|\frac{\mu_n(t)}{f(n)}-\mu(\lambda_n(t))\right|< \varepsilon,\ \
\sup_{0\le t,s\le 1, s\not=t}\left|ln\frac{\lambda(t)-\lambda(s)}{t-s}\right|< \varepsilon.$$
Then for all $t\in[0,1]$ and $n>n_2$  holds the inequality: $\mu_n(t)>(c-\varepsilon)f(n).$
Thus $f(n)\to\infty$ as $n\to\infty$, it exists $n_3$ such that for all $n>n_3$ we can assume $f(n)>\frac{\max\{n_1,n_2\}}{c-\varepsilon}$.

Then for all $n>n_3$ 
 $$\rho_D(x_{\mu_n},x)<\varepsilon,$$
which implies the lemma assertion.

{\it Proof of the theorem.}
By Skorokhod theorem about one probability space (see, for example, theorem 11 in section V in \cite{Bulinski}), it exists the probability space $(\Omega'_1,{\mbox{\gothic{A}}}_1 , P_1)$ and the random processes $X_n:\Omega'_1\rightarrow D[0,1]$ such that

\medskip

$1) X_n\stackrel{d}{=}Y_n, \ \ X\stackrel{d}{=}Y;$

\medskip

$ 2) X_n\stackrel{a.s.}{\to} X \mbox{ at } n\to\infty \mbox{ in }
D[0,1]. $

\medskip

Denote by $\Omega_1$ the measurable subset of
$\Omega'_1$ such that $P_1(\Omega_1)=1$ and the convergence 2) is true for all $\omega_1\in\Omega_1$.

 Also by Skorokhod theorem the convergence $\frac{\nu_n}{f(n)}\stackrel{d}{\to}\nu$ at $n\to\infty$ in
$D[0,1]$ implies the existence of probability space
$(\Omega'_2, {\mbox{\gothic{A}}}_2, P_2)$
 and random processes $\mu_n:\Omega'_2\rightarrow D[0,1]$ such that

\medskip

 $ 1) \nu_n\stackrel{d}{=}\mu_n, \ \ \nu\stackrel{d}{=}\mu;$

\medskip

 $ 2) \frac{\mu_n}{f(n)}\stackrel{a.s.}{\to} \mu \mbox{ as } n\to\infty \mbox{ in } D[0,1]. $

\medskip

Note that the set $$A=\{f\in D[0,1] : f(t)\ge 0, \ f(t) \mbox{ is a non-decreasing  function, taking values in }{\rm I\kern-.23em N}\}$$ is measurable relating Borel $\sigma$ - algebra in D[0,1]. Then $1=P(\nu_n\in A)=P(\mu_n\in A)$, that is almost all  trajectories of $\mu_n$ are also non-negative, non-decreasing and take values in {\rm I\kern-.23em N}. Similarly we can see that almost sure $\mu(t)>c$ for all $t \in [0,1]$.

Denote by $\Omega_2$ the measurable subset of $\Omega'_2$ such that $P_2(\Omega_2)=1$, all trajectories of $\mu_n$ are non-negative and non-decreasing,  the convergence  2) and inequality  $\mu(t)>c$, $t \in [0,1]$ are true for all  $\omega_2\in\Omega_2$.

 Further we will consider the probability space
 $(\Omega,  {\mbox{\gothic{A}}}, P)$, where $\Omega=\Omega_1\times\Omega_2$, the $\sigma$-algebra $ {\mbox{\gothic{A}}}$ consists of elements of $\sigma$-algebra $ {\mbox{\gothic{A}}}_1\times {\mbox{\gothic{A}}}_2$
belonging to $\Omega$, the probability $P$ is a restriction of probability $P_1\otimes P_2$ to $\sigma$-algebra $ {\mbox{\gothic{A}}}$.

 Let $\omega=(\omega_1,\omega_2)\in \Omega$. Consider the sequence $x_n\equiv X_n(\omega_1),$ $ x\equiv X(\omega_1),$    $ \gamma_n(t)\equiv\mu_n(t)(\omega_2), $   $\gamma(t)\equiv \mu(t)(\omega_2)$. 

Then by Lemma \ref{lemma} we obtain the convergence $x_{\gamma_n}\to x$ as $n \to\infty$ in $D[0,1]$ for all $\omega\in\Omega$,  which implies the convergence $X_{\mu_n}\stackrel{a.s.}{\to}X$ as $n \to\infty$ in $D[0,1]$. Because the distributions of $X_{\mu_n}$ and $Y_{\nu_n}$ coincide, the assertion of theorem is proved.

\begin{rem}
Note, that Theorem \ref{maintheorem} not holds if trajectories of $Y$ not continuous.
\end{rem}

\begin{exm}
Consider the sequence of non-random functions
 $$
 x_{2n}(t)=\left\{
 \begin{array}{lr} 0
\mbox{ if }&  0\le t\le \frac{1}{2}-\frac{1}{2^{2n}},\\
2^{2n} t+1-2^{2n-1}\mbox{ if }& \frac{1}{2}-\frac{1}{2^{2n}}< t\le\frac{1}{2},\\
1\mbox{ if }& \frac{1}{2}< t\le 1,
\end{array}
\right.
$$
$$
 x_{2n+1}(t)=\left\{
 \begin{array}{lr} 0
\mbox{ if }&  0\le t\le \frac{1}{2},\\
2^{2n+1} t-2^{2n}\mbox{ if }& \frac{1}{2}< t\le\frac{1}{2}+\frac{1}{2^{2n+1}},\\
1\mbox{ if }& \frac{1}{2}+\frac{1}{2^{2n+1}}< t\le 1.
\end{array}
\right.
$$
We will define $\lambda_n$ as:

$$
 \lambda_{2n}(t)=\left\{
 \begin{array}{lr}
(1-\frac{1}{2^{2n-1}})t\mbox{ if }& 0\le t\le\frac{1}{2},\\
(1+\frac{1}{2^{2n-1}})t- \frac{1}{2^{2n-1}}\mbox{ if }& \frac{1}{2}< t\le 1.
\end{array}
\right.
$$
$$
 \lambda_{2n+1}(t)=\left\{
 \begin{array}{lr}
(1+\frac{1}{2^{2n}})t\mbox{ if }& 0\le t\le\frac{1}{2},\\
(1-\frac{1}{2^{2n}})t+ \frac{1}{2^{2n-1}}\mbox{ if }& \frac{1}{2}< t\le 1.
\end{array}
\right.
$$
It is easy to see that $$\sup_{0\le t\le 1}|x_n(\lambda_n(t))-x(t)|\to 0, $$
where $$x(t)=\left\{
 \begin{array}{lr}
0\mbox{ if }& 0\le t <\frac{1}{2},\\
1 \mbox{ if }& \frac{1}{2}\le t\le 1.
\end{array}
\right.
$$
Let the sequence of random variables 
$$\nu_n=\left\{
 \begin{array}{lr}
2n\mbox{ with probability }& \frac{1}{2},\\
2n+1 \mbox{ with probability }& \frac{1}{2}.
\end{array}
\right.
$$
Then for arbitrary $\mu_n\in \Delta[0,1]$ it holds:
\begin{equation}  \label{example}
\sup_{0\le t\le 1}|x_{\nu_n}(\mu_n(t))-x(t)|=\frac{1}{2}\left(\sup_{0\le t\le 1}|x_{2n}(\mu_n(t))-x(t)|+ \right.
\end{equation}
$$\left.\sup_{0\le t\le 1}|x_{2n+1}(\mu_n(t))-x(t)|\right).$$
We will show there is no such $\mu_n\in\Delta[0,1]$ that (\ref{example}) tands to zero.

It is easy to see that $\mu_n(\frac{1}{2})$ can't be equal to $\frac{1}{2}$. Assume that $\mu_n(\frac{1}{2})<\frac{1}{2}$. Then $x_{2n+1}(\mu_n(\frac{1}{2}))=0$ and $x(\frac{1}{2})=1$, which  makes convergence impossible. Similarly, if $\mu_n(\frac{1}{2})>\frac{1}{2}$ then $x_{2n}(\mu_n(\frac{1}{2}))=0$.
\end{exm}

In the following we will consider some corollaries of main result.

\begin{cor}
Let $X_i$ be the independent random variables such that $\sum_{i=1}^n X_i\stackrel{d}{\to}X$ and $\pi(t)$ is an independent of $X_i \ (i\in {\rm I\kern-.23em N})$ Poisson random process, $E\pi(t)=t$ and  $a>0$. Then it holds
$$\sum_{i=1}^{\pi(n(t+a))}X_i\stackrel{d}{\to}X \mbox{ as } n\to \infty \mbox{ in } D[0,1]. $$
\end{cor}
{\it Proof.} Note that $\frac{\pi(n(t+a))}{n}\stackrel{d}{\to} t+a \mbox{ in } D[0,1] \mbox{ as } n\to\infty. $ Thus  all conditions of  Theorem \ref{maintheorem} are satisfied.

Let $X_i$ be the i.i.d. random variables such that $E X_i=0, D X_i=\sigma^2$  and $$S_i=\sum_{k=1}^i X_k. $$ 
Consider the random processes
$$X_n(t)=\frac{1}{\sigma \sqrt{n}}S_{i-1}+\frac{t-(i-1)/n}{1/n}\frac{1}{\sigma \sqrt{n}}X_i, \ t\in\left[\frac{i-1}{n},\frac{i}{n}\right]. $$ It is well-know that $X_n\stackrel{d}{\to}W $ as $n\to\infty$ in $D[0,1]$ (here by $W$ denoted a Wiener random process).

\begin{cor} $[$generalized invariance principle$]$
Let $\nu_n$ is a sequence of random elements 	
taking values in  ${\rm I\kern-.23em N}$, $f(n)$ is a function taking values in ${\rm I\kern-.23em R}$ such that $f(n)\to\infty$ as $n\to\infty$ and
$$ \frac{\nu_n}{f(n)}\stackrel{d}{\to} \nu$$ and it exists $c>0$ such that $\nu>c. $ Then it holds $$X_{\nu_n}\stackrel{d}{\to}W \mbox{ as } n\to\infty \mbox{ in } D[0,1], $$
where $$X_{\nu_n}(t)=\frac{1}{\sigma \sqrt{\nu_n}}S_{i-1}+\frac{t-(i-1)/\nu_n}{1/\nu_n}\frac{1}{\sigma \sqrt{\nu_n}}X_i, \ t\in\left[\frac{i-1}{\nu_n},\frac{i}{\nu_n}\right]. $$
\end{cor}

\begin{rem}
Random process $X (t)$ is usually interpreted as a process of random walk of a particle, which changes direction at moments of times $t=\frac{i}{n}$. Random process $X_{\nu_n}$ can be interpreted as a random walk of a particle, which changes direction at random moments of time.
\end{rem}

\section{Estimation of convergence rate}
\begin{thm}
Let $Y_n, Y$ and $\nu_n, \nu$ are the  random processes in $D[0,1]$, $\nu_n(t)$ has non-decreasing non-negative trajectories, for every $t\in [0,1]$ $\nu_n(t)$  	
take values in  ${\rm I\kern-.23em N}$, $f(n)$ is a function taking values in ${\rm I\kern-.23em R}$  and besides it exists $c>0$ such that  $\nu(t)>c$ for all $t \in [0,1]$.  Let $Y_n$ and $\nu_n$ are independent.
Then
$$\sup_{x}\left|P\{Y_{\nu_n}<x\}-P\{Y<x\}\right| \le \sup_x\sup_{k\ge [cf(n)]}\left|P\{Y_k<x\}-P\{Y<x\}\right|+$$ $$+2\sup_x\sup_k\left|P\{Y_k<x\}-P\{Y<x\}\right|\inf_x\left|P\{\frac{\nu_n}{f(n)}<x\}-P\{\nu<x\}\right|.$$
\end{thm}
{\bf Proof.} 
Let $x\in {\rm I\kern-.23em R}$ is arbitrary. Then
$$
\left|P\{Y_{\nu_n}<x\}-P\{Y<x\}\right|=\left|\sum_{k=1}^{\infty}(P\{Y_k<x\}-P\{Y<x\})P\{\nu_n=k\} \right|\le$$
$$\sum_{k=1}^{\infty}|P\{Y_k<x\}-P\{Y<x\}|P\{\nu_n=k\}\le$$
\begin{equation}                \label{ineq}
\sum_{k=1}^{\infty}|P\{Y_k<x\}-P\{Y<x\}|P\left\{\ \frac{k-1}{f(n)}<\nu\le \frac{k}{f(n)}\right\}+ \end{equation}

$$\left|\sum_{k=1}^{\infty}|P\{Y_k<x\}-P\{Y<x\}|\left(P\left\{\frac{k-1}{f(n)}<\nu_n\le \frac{k}{f(n)}\right\}-P\left\{ \frac{k-1}{f(n)}<\nu\le \frac{k}{f(n)}\right\}\right)\right|.$$ 

Note that for all  $k<cf(n)$  the equality $P\{\nu\le\frac{k}{f(n)}\}=0$ is holds. Then the first term of (\ref{ineq}) can be estimated by
$$\sum_{k=1}^{\infty}|P\{Y_k<x\}-P\{Y<x\}|P\left\{\ \frac{k-1}{f(n)}<\nu\le \frac{k}{f(n)}\right\}\le\max_{k\ge cf(n)}|P\{Y_k<x\}-P\{Y<x\}|. $$
Let $N\in {\rm I\kern-.23em N}$ is arbitrary. For the second term of (\ref{ineq}) we have
$$\left|\sum_{k=1}^{\infty}|P\{Y_k<x\}-P\{Y<x\}|\left(P\left\{\frac{k-1}{f(n)}<\nu_n\le \frac{k}{f(n)}\right\}-P\left\{ \frac{k-1}{f(n)}<\nu\le \frac{k}{f(n)}\right\}\right)\right|$$
$$\le\sup_{k}|P\{Y_k<x\}-P\{Y<x\}|\left|P\left\{\nu_n\le \frac{N}{f(n)}\right\}-P\left\{ \nu\le \frac{N}{f(n)}\right\}\right|+ $$
$$\sup_{k}|P\{Y_k<x\}-P\{Y<x\}|\left|P\left\{\nu_n> \frac{N}{f(n)}\right\}-P\left\{ \nu> \frac{N}{f(n)}\right\}\right|=  $$
$$=2 \sup_{k}|P\{Y_k<x\}-P\{Y<x\}|\left|P\left\{\nu_n\le \frac{N}{f(n)}\right\}-P\left\{ \nu\le \frac{N}{f(n)}\right\}\right|.  $$
The arbitrariness of $N$ implies the following inequality for  second term of (\ref{ineq})
$$\left|\sum_{k=1}^{\infty}|P\{Y_k<x\}-P\{Y<x\}|\left(P\left\{\frac{k-1}{f(n)}<\nu_n\le \frac{k}{f(n)}\right\}-P\left\{ \frac{k-1}{f(n)}<\nu\le \frac{k}{f(n)}\right\}\right)\right| $$
$$\le  2\sup_k\left|P\{Y_k<x\}-P\{Y<x\}\right|\inf_x\left|P\{\frac{\nu_n}{f(n)}<x\}-P\{\nu<x\}\right|.$$
That completes  the proof.


\end{document}